\documentclass{article}
\usepackage{epsfig}
\title{Slopes of a Hypergeometric System \\
Associated to a Monomial Curve}
\author{Francisco Jes\'{u}s Castro-Jim\'{e}nez\footnote{Partially supported by DGESIC and FQM-218} and Nobuki Takayama}
\date{July 4, 2001}
\def\pd#1{ \partial_{#1} }
\newenvironment{proof}{\par\noindent Proof:\ }{${\tt [}\kern-0.2mm{\tt ]}$}
\newtheorem{theorem}{Theorem}[section]
\newtheorem{proposition}[theorem]{Proposition}
\newtheorem{definition}[theorem]{Definition}
\newtheorem{example}[theorem]{Example}
\newtheorem{algorithm}[theorem]{Algorithm}

\newcommand{\CC}{{\bf C}}
\newcommand{\cD}{{\cal D}}
\newcommand{\cH}{{\cal H}}
\newcommand{\cO}{{\cal O}}

\newcommand{\RR}{{\bf R}}
\newcommand{\gr}{{\rm gr}}
\newcommand{\inn}{{\rm in}}

\begin{document}
\maketitle

\section{Introduction}
Let $A_n={\bf C}\langle x_1,\ldots,x_n,\partial_1,\ldots,\partial_n\rangle$
be the Weyl algebra of order $n$ over the complex numbers $\bf C$ and
${\bf C}[\partial]={\bf C}[\partial_1,\ldots,\partial_n]$
the subring of $A_n$ of linear differential operators with constant
coefficients.

Let $A=(a_{ij})$ be an integer $d\times n$-matrix of rank $d$.
We denote by $I_A \subset {\bf C}[\partial]$ the {\it toric ideal}
associated to $A$: {\it i.e.} $I_A$ is the ideal generated by the
set $$ \{\partial^u-\partial^v\, \vert \, u,v\in {\bf N}^n,\,
Au^T=Av^T\}$$ where $()^T$ means ``transpose".

We denote by $\theta$ the vector $(\theta_1,\ldots,\theta_n)^T$ with
$\theta_i=x_i\partial_i$. For a given $\beta =
(\beta_1,\ldots,\beta_d)^T \in {\bf C}^d$ we consider the column
vector (in $A_n^d$) $A\theta-\beta$  and we denote by $\langle
A\theta-\beta\rangle$ the left ideal of $A_n$ generated by the
entries of $A\theta-\beta$.

Following Gel'fand, Kapranov and Zelevinsky \cite{GKZ},
we denote by $H_A(\beta)$
the left ideal of $A_n$ generated by $I_A \cup \langle
A\theta-\beta \rangle$.
It is called the GKZ-hypergeometric
system associated to the pair $(A,\beta)$.
The quotient  ${\cal
H}_A(\beta)=A_n/H_A(\beta)$ is a holonomic $A_n$-module (see e.g.
\cite{SST}).

If the toric ideal $I_A$ is homogeneous, i.e., if  the
 ${\bf Q}$-row span of $A$  contains $(1,\ldots,1)$, it is
known (\cite{Hotta}, see also \cite{SST}) that ${\cal H}_A(\beta)$
is regular holonomic and the book \cite{SST} is devoted to
an algorithmic study of such systems.
Especially, the book gives an algorithmic method to construct
series solutions around singular points of the system.

In this article, we would start a study of singularities
of GKZ-hypergeometric systems
for non homogeneous toric ideals $I_A$ by treating the ``first"
case when $d=1$, $A=(a_1, a_2, \ldots, a_n)\in {\bf Z}^n$, $a_1=1$.
We evaluate the geometric slopes of ${\cal H}_A(\beta)$ by successive
restrictions of the number of variables.
The slopes characterize Gevrey class solutions
around singular locus.

Our evaluation of the slopes is done as follows:
(1) We translate Laurent and Mebkhout's theorem \cite{Lau-Meb-2}
on restrictions and slopes of ${\cal D}$-modules
into an algorithm to evaluate the slopes
by utilizing the results \cite{ACG} and \cite{O}.
(2) Apply our general algorithm to the hypergeometric system associated
to $A=(1,a_2, \ldots, a_n)$.
This system has a lot of nice properties and our algorithm outputs the slopes
without computation on computers.

\section{Micro-characteristic varieties}
\label{section:microcharacteristic} \setcounter{equation}{0}

In this section following Laurent \cite{Lau-ens}, we describe
micro-characteristic varieties for a given $\cD$-module. We will
state a result of Laurent and Mebkhout \cite[Corollaire 2.2.9
]{Lau-Meb-2} (see also (\cite{Meb-positivite}, p. 125) and
(\cite{Lau-Meb-papa}, p. 42)), allowing
to reduce our general problem of evaluating the slopes to less variables.

In this section $X={\bf C}^n$ and ${\cal D}_X={\cal D}$ is the
sheaf of linear differential operators with holomorphic function
coefficients. Let $M$ be a coherent $\cD_X$-module. Recall that
the characteristic variety of $M$ (denoted by ${\rm Ch}(M)$) is a
analytic sub-variety of the cotangent bundle $T^* X$.

Suppose that $Y\subset X$ is a smooth hypersurface.
We say that $Y$ is {\it non-characteristic} for $M$ if $T^*_Y
X \cap {\rm Ch}(M) \subset T^*_XX$.
Here $T^*_Y X$ is the conormal
bundle to $Y$ in $X$ and $T^*_XX$ is the zero section of
$T^*X$.

Now, following Laurent \cite{Lau-these,Lau-ens}, we shall
define the notion of {\it non-micro-characteristic variety} for $M$.
To simplify the presentation we will assume that
$(x_1,x_2,\ldots,x_n)$ are local coordinates in $X$ and  $Y$
is defined by $x_n=0$. We denote by
$(x_1,\ldots,x_n,\xi_1,\ldots,\xi_n)$ local coordinates in
$T^*X$. Sometimes it will be useful to write $x_1=y_1, \ldots,
x_{n-1}=y_{n-1}$, $x_n=t$, $\xi_1=\eta_1, \ldots,
\xi_{n-1}=\eta_{n-1}$ and $\xi_n=\tau$.

Let us denote by $\Lambda$
the conormal bundle of $Y$ in $T^*X$ (i.e. $\Lambda =
T^*_YX$). So, in local coordinates $\Lambda = \{(y,t,\eta,\tau)\in
T^*X\, | \, t=\eta=0\}$.
Here,
$\eta=(\eta_1,\ldots,\eta_{n-1})$ and $y=(y_1,\ldots,y_{n-1})$. We
denote by $(y,\tau,y^*,\tau^*)$ local coordinates on the cotangent
bundle $T^*\Lambda$.

We denote by $V_\bullet(\cD)$ (or simply by $V$) the
Malgrange-Kashiwara filtration associated to $Y$ on $\cD$ and by
$F_\bullet(\cD)$ (or simply by $F$) the order filtration on $\cal
D$. For a given rational number $p/q\geq 0$ we denote by $L_{p/q}$
the filtration on $\cD$ defined by $pF+qV$.
The $L_{p/q}$-order of a monomial $y^\alpha
t^l\partial_y^\beta\partial_t^k$ is equal to
$p(|\beta|+k)+q(k-l)$, where
$|\beta|=\beta_1+\cdots+\beta_{n-1}$.
We will simply write $L=L_{p/q}$ if no confusion arises.

For $p>0$ the associated graded ring $\gr^{L}(\cD)$ is canonically
isomorphic to $\pi_*\cO_{[T^*\Lambda]}$ \cite[p. 407]{Lau-ens},
where $\pi: T^*\Lambda \rightarrow \Lambda$ is the canonical projection
and $\cO_{[T^*\Lambda]}$ denotes holomorphic functions on
$T^*\Lambda$ which are polynomials on the fibers of $\pi$.
In local coordinates,
$\gr^L(\cD_0)$ is expressed as
${\bf C}\{y\}[\tau,y^*,\tau^*]$ where $\cD_0$ is the stalk of $\cD$ at
the origin.

Given a differential operator
$$P=\sum_{\alpha l\beta k} p_{\alpha l\beta k}y^\alpha t^l\partial_y^\beta\partial_t^k$$
the $L$-order
of $P$ is the maximun value of $p(|\beta|+k)+q(k-l)$ over the
monomials of $P$.
For $p>0$ we define the $L$-principal symbol of
$P$ by $$\sigma^{L}(P)= \sum p_{\alpha l\beta k}y^\alpha
(\tau^*)^l(y^*)^\beta(-\tau)^k$$ where the sum is taken over
monomials with maximal $L$-order. The $L$-principal symbol of $P$
is an element of $\gr^{L}(\cD)$ and then is a function on
$T^*\Lambda$.
In the classical case, i.e., for $L=F$,
$\gr^F(\cal D)$ is identified with
${\bf C}\{x\}[\xi_1,\ldots,\xi_n] =
 {\bf C}\{y,t\}[\eta,\tau]$
and the $F$-principal symbol of $P$ is
simply denoted by $\sigma^F(P)=\sum p_{\alpha l\beta k}y^\alpha
t^l \eta^\beta \tau^k$ where the sum is taken for $|\beta|+k$
maximum.

For each left ideal $I\subset \cal D$ we denote by $\sigma^L(I)$ the
ideal of $\gr^L(\cal D)$ generated by the set of $\sigma^L(P)$ for
$P\in I$.

For each $L$-filtration on $\cD$ we associate a ``good"
$L$-filtration on $M$, by means of a finite presentation. The
associated $\gr^{L}(\cD)$-module $\gr^{L}(M)$ is coherent (see
\cite[3.2.2]{Lau-ens}). The radical of the annihilating ideal
${\rm Ann}_{\gr^{L}(\cD_X)}(\gr^{L}(M))$,
which is independent of the ``good" filtration on $M$,
defines an analytic sub-variety of $T^*\Lambda$.
This variety is called the
$L$-{\it characteristic variety} of $M$ and it is denoted by
${\rm Ch}^{L}(M)$.

Suppose now $Z\subset X$ is a smooth hypersurface transverse to
$Y$.
Suppose for simplicity $Z$ is defined in local coordinates by
$y_1=0$. The conormal space $\Lambda' := T^*_{Y\cap Z}Z$ is a
smooth subvariety of $\Lambda= T^*_Y X$ defined in local
coordinates by $y_1=0$. So $T^*_{\Lambda'}\Lambda$ is the
subvariety of $T^*\Lambda$ defined in local coordinates by
$y_1=y'^{*}=\tau^*=0$, where $y'=(y_2,\ldots,y_{n-1})$.

\begin{definition} {\rm \cite{Lau-these}, \cite{Lau-Meb-2}}
\rm
We say that $Z$ is {\it non-micro-characteristic of type} $L$ for
$M$ if $T^*_{\Lambda'}\Lambda \cap {\rm Ch}^{L}(M)$ is contained in
$T^*_\Lambda \Lambda$.
Sometimes we will say that, if this
condition holds, $Z$ is {\it non $L$-micro-characteristic} for $M$.
\end{definition}

The sheaf of rings $\gr^L(\cD)$ is endowed with two graduations:
First one is induced by the $F$-filtration  and the second one is
induced by $V$-filtration. Recall Laurent's definition of slope of
a coherent $\cal D$-module $M$.

\begin{definition}{\rm \cite{Lau-ens}}\label{laurent-slope} \rm
The rational number $-p/q$ is said to be a {\it slope} of $M$ w.r.t.
$Y$ at the origin if and only if the radical of the ideal
${\rm Ann}_{\gr^{L}(\cD_X)}(\gr^{L}(M))$ is not bi-homogeneous
for $F$ nor $V$ graduations.
\end{definition}

\noindent
An important consequence of the work \cite{Lau-Meb-papa} (see Th. 2.4.2)
is what  follows:  An holonomic ${\cal D}_X$-module $M$  is regular w.r.t.
$Y$ at the origin if and only if $M$ has no slope w.r.t. $Y$ at the origin. We
will use this fact freely in the text.

\noindent
{\it Remark}:
In \cite{Lau-Meb-2}, $\infty$ and $0$ ($F$ and $V$) are
included in the set of the slopes.
We do not include them in the set of the slopes in this paper.
\smallbreak

Finally the following result by Laurent and Mebkhout
allows induction on the number of variables to calculate slopes.

\begin{theorem}{\rm \cite[Corollaire 2.2.9]{Lau-Meb-2}}
{\rm  \label{Cauchy-Kow}}
Let $M$ be a holonomic $\cD_X$-module. Let $Z$ and $Y$  be
transverse smooth hypersurfaces on $X$ such that $Z$  is non
$L$-micro-characteristic {\rm (}for all $p>0, q > 0${\rm )} for $M$.
Then the slopes of $M$ w.r.t. $Y$  equals the slopes of $M'$ w.r.t.
$Y'=Z\cap Y$ where $M'$ is the restriction of $M$ to $Z$.
\end{theorem}

This is a deep result in $\cal D$-module theory. Its proof uses
the algebraic-analytic comparison theorem (\cite[Th. 2.4.2
]{Lau-Meb-papa}) and a Cauchy-Kowalewska theorem for Gevrey
functions w.r.t. $Z$ (\cite[Corollaire 2.2.4]{Lau-Meb-2}; see also
\cite[Th. 6.3.4]{Meb-positivite}.)

\section{Computing slopes by reducing the number of variables}
\label{section:general-algorithm}
\setcounter{equation}{0}

We have introduced the notion of the slopes and the invariance of
them under restrictions satisfying a condition on $L$-characteristic
varieties.
In the sequel, we assume that
our ideal is that of the Weyl algebra $A_n$.  Constructions in
sheaves such as restrictions, micro-characteristic varieties in
the previous section can be done via constructions in the Weyl
algebra as we usually see in the computational $\cD$-module
theory.

We are interested in computation of the slopes. The slopes of
$A_n/I$ (at the origin) along $x_n=0$ can be computed by the ACG
algorithm introduced in \cite{ACG}.
In this section, translating Laurent and Mebkhout's result into
computer algebra algorithms, we will give
a preprocessing method for the ACG algorithm to accelerate the
original. The preprocessing is useful for a class of inputs
including GKZ hypergeometric ideals as we will see in Section
\ref{section:computingSlopes}.
Let us firstly recall the ACG algorithm.

A weight vector is an element
$W=(u_1,\ldots,u_n,v_1,\ldots,v_n)\in \RR^{2n}$ such that
$u_i+v_i\geq 0$ for all $i$.
This weight vector $W=(u,v)$
induces a natural filtration on $A_n$ and it is called the
$W$-filtration.
The associated graded ring is denoted by
$\gr^W(A_n)$ and for each left ideal $I\subset A_n$ the
associated graded ideal is denoted by $\inn_W(I)$ or
$\inn_{(u,v)}(I)$.
Here, the {\it initial ideal} $\inn_{(u,v)}(I)$ is the ideal
generated by $\inn_{(u,v)}(f), f \in I$ in ${\rm gr}^W(A_n)$.
When $u_i + v_i >0$, it is an ideal in the polynomial ring
of $2n$ variables:
${\rm gr}^W(A_n) = {\bf C}[x_1, \ldots, x_n, \xi_1, \ldots, \xi_n]$.
The initial ideal of $I$ with respect to the weight $(u,v)$
is generated by the $(u,v)$-initial terms
of a Gr\"obner basis of $I$ by an order which refines the partial
order defined by $(u,v)$.
See, e.g., \cite[Theorem 1.1.6]{SST}.

Consider the filtration $L = p F + q V, p > 0, q > 0$
introduced in the previous section.
The ideal $\sigma^L(I)$, which gives the $L$-characteristic variety,
can be expressed in terms of the initial ideal as follows:
$$ \sigma^L(I) = {\rm gr}^L({\cal D}) \cdot \inn_\ell(I)
  |_{x_1 \mapsto y_1, \ldots, x_{n-1} \mapsto y_{n-1},x_n \mapsto \tau^*,
   \xi_1 \mapsto y_1^*, \ldots, \xi_{n-1} \mapsto y_{n-1}^*,\xi_n\mapsto-\tau},
$$
$ \ell = p ( \overbrace{0, \ldots, 0}^{n}, \overbrace{1, \ldots, 1}^{n})+
         q ( \overbrace{0, \ldots, 0, -1}^{n}, \overbrace{0, \ldots,0,1}^{n})
$.

For two weight vectors $W$ and $W'$ and a term order $<$, we
denote by $<_{W,W'}$ the order
\begin{eqnarray*}
&& x^\alpha \pd{}^\beta <_{W,W'} x^a \pd{}^\beta\\
&\Leftrightarrow& W \cdot (\alpha,\beta) < W \cdot (a,b) \\ &
& \mbox{or}\, W \cdot (\alpha,\beta) = W \cdot (a,b) \
      \mbox{ and } W' \cdot (\alpha,\beta) < W' \cdot (a,b) \\
&&  \mbox{or}\, W^* \cdot(\alpha,\beta) = W^* \cdot (a,b)\
\mbox{ for both }\
     W^*=W, W' \ \mbox{ and} \
    x^\alpha \pd{}^\beta < x^a \pd{}^b.
\end{eqnarray*}

To each differential operator $P=\sum
p_{\alpha\beta}x^\alpha\partial^\beta\in A_n$ we associate the
Newton polygon $N(P)$ of $P$ (w.r.t. $x_n=0$) defined as the
convex hull of the subset of ${\bf Z}^2$
$$\bigcup_{p_{\alpha\beta}\not= 0}(|\beta|,\beta_n-\alpha_n)+
(-{\bf N})^2.$$

Let $I$ be a left ideal in $A_n$. As we said in
Definition \ref{laurent-slope} the notion of slope of a differential
system was introduced by Y. Laurent \cite{Lau-ens}.
Let us give here a slight different but equivalent definition:
the number $r$,
$-\infty < r < 0$ is a {\it geometric} slope of $I$
(or of $A_n/I$) w.r.t. $x_n=0$ if and only if
$\sqrt{\sigma^{(-r)F+V}(I)}$ is not bihomogeneous with respect to the
weight vectors $F=(0,\ldots,0,1,\ldots,1)$ and
$V=(0,\ldots,0,-1,0,\ldots,0,1)$.
Following \cite{ACG}, we say
that the number $r$, $-\infty < r < 0$ is an {\it
algebraic} slope of $I$ (or of $A_n/I$) if and only if
$\sigma^{(-r)F+V}(I)$ is not bihomogeneous with respect to the
weight vectors $F$ and $V$.
The geometric slope is simply called the slope in this paper
if confusion does not arise.
Note that we may consider $\inn_L(I)$ instead of $\sigma^L(I)$
as far as we are concerning about homogeneity.
For algebraic or geometric slope $r$, the
weight vector $L=(-r)F+V$ lies on a face of the Gr\"obner fan
of $I$ (\cite{ACGjpaa}, \cite{SST}), which yields the
following algorithm.

\begin{algorithm} \rm (\cite{ACG}, ACG algorithm) \\
Input: $G = \{ P_1, \ldots, P_m \}$ (generators of an ideal $I$) \\
Output: All algebraic and geometric slopes of $A_n/I$ w.r.t. $x_n=0$ at the
origin.\\
${\tt geometric\_slope} = \emptyset$; ${\tt algebraic\_slope} = \emptyset$; \\
${\tt F}=(0, \ldots, 0, 1,\ldots, 1)$; ${\tt V} = (0, \ldots, 0, -1,0, \ldots, 0, 1)$;\\
${\tt p}=1$; ${\tt q}=0$; ${\tt slope}=-\infty$;
${\tt previous\_slope} = {\tt slope}$; \\
{\tt while} $({\tt slope}!= 0)$  {\tt \{} \\
  \mbox{\,}\quad ${\tt L} = {\tt p F} + {\tt q V}$; \\
  \mbox{\,}\quad ${\tt G} = \mbox{a Gr\"obner basis of } I
                          \mbox{ with respect to the order } <_{L,V}$; \\
  \mbox{\,}\quad ${\tt slope} = {\rm the\  minimum}\  \mbox{ of } \ 0
                                 \ \mbox{ and } \\
  \mbox{\,}\quad \quad\quad \{\mbox{ the slopes $r$ of the Newton polygon }
          N(P) \,|\, P \in G, r > {\tt previous\_slope} \}$ \\
  \mbox{\,}\quad {\tt if} ${\tt slope} = 0$ ,
        then {\tt return}({\tt algebraic\_slope} and {\tt geometric\_slope}). \\
  \mbox{\,}\quad {\tt if}
    $\sigma^{L}(G)$ is not homogeneous for $F$ nor $V$ {\tt then} {\tt \{} \\
        \mbox{\,}\quad\quad ${\tt algebraic\_slope}
         = {\tt algebraic\_slope} \cup \{ {\tt slope} \}$  \\
        \mbox{\,} \quad {\tt \}} \\
   \mbox{\,}\quad {\tt if} $\sqrt{\sigma^{L}(G)}$ is not homogeneous
             for $F$ nor $V$ {\tt then} {\tt \{} \\
          \mbox{\,}\quad\quad ${\tt geometric\_slope}
         = {\tt geometric\_slope} \cup \{ {\tt slope} \}$  \\
             \mbox{\,}\quad {\tt \}} \\
   \mbox{\,}\quad ${\tt p} = {\rm numerator}(|{\tt slope}|);
   {\tt q} = {\rm denominator}(|{\tt slope}|);$  \\
   \mbox{\,}\quad ${\tt previous\_slope} = {\tt slope};$ \\
  {\tt \}} \\
  {\tt return} ({\tt algebraic\_slope} and {\tt geometric\_slope}) \\
\end{algorithm}

\noindent
Here, we use the convention
$F:=\infty F + V$.
Note that the ideal $J = \langle f_1, \ldots, f_m \rangle
\subset {\bf C}[x_1, \ldots, x_n, \xi_1, \ldots, \xi_n]$
is homogeneous for the weight $(u,v) \in {\bf Z}^{2n}$
if and only if all $(u,v)$-homogeneous subsums of $f_i$ belong
to the ideal $J$.

\bigbreak

\noindent
For ordinary differential equations, the ACG algorithm is nothing
but the well-known
Newton polygon method; see two examples below.
\begin{example} \rm
Put $G = \{ x_1^{p+1} \pd{1} + p \}$, $n=1$, $p \in {\bf N}$.
Then, the ACG algorithm returns $\{-p\}$.
( $\exp(x_1^{-p})$ is a classical solution of $G$.)
\end{example}
\begin{example} \rm
Put $G = \{ 2 x_1 (x_1 \pd{1})^2 + x_1^2 \pd{1} +1 \}$, $n=1$.
Then, the ACG algorithm returns $\{-1/2\}$.
( $\exp(x_1^{-1/2})$ is a classical solution of $G$.)
\end{example}

\noindent
Next example of two variables is generated by a computer algebra system
{\tt kan/k0} \cite{openxm}.
\begin{example} \rm
We consider
the GKZ hypergeometric ideal $I$ associated to the matrix
$A=(1,3)$ and $\beta=-3$.
We will compute the slopes of $A_2/I$ at the origin along $x_2=0$
by the ACG algorithm.
The Gr\"obner basis of $I$
for the weight
$(\overbrace{0}^{x_1},\overbrace{0}^{x_2},
  \overbrace{1}^{\partial_1},\overbrace{1}^{\partial_2}) $ is
{\footnotesize
\begin{verbatim}
[  x1*Dx1+3*x2*Dx2+3 , -Dx1^3+Dx2 , -3*x2*Dx1^2*Dx2-5*Dx1^2-x1*Dx2 ,
 -9*x2^2*Dx1*Dx2^2-36*x2*Dx1*Dx2+x1^2*Dx2-20*Dx1 ,
 -27*x2^3*Dx2^3-189*x2^2*Dx2^2-x1^3*Dx2-276*x2*Dx2-60 ]
\end{verbatim} }
Here, {\tt Dx}$i$ and  {\tt x}$i$ stand for $\partial_i$ and $x_i$
respectively.
The Newton polygons $N(P)$'s are figured below.

\begin{center}
\epsfxsize=10cm \epsffile{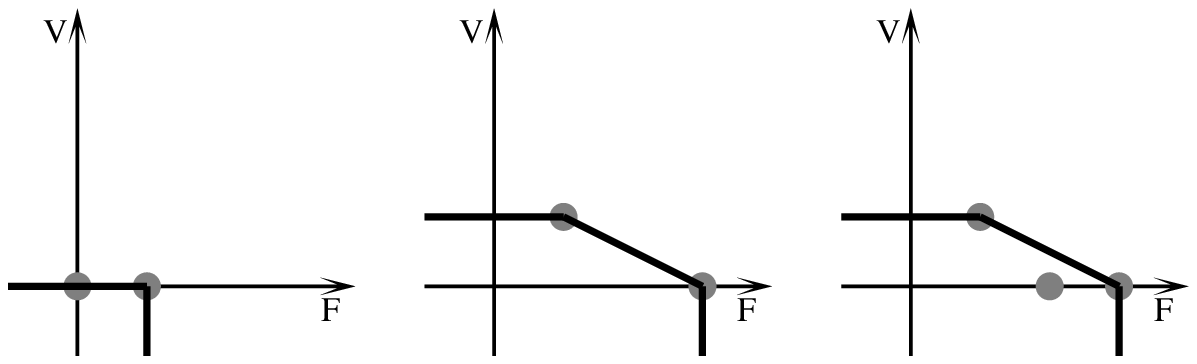}
\end{center}
\begin{center}
\epsfxsize=10cm \epsffile{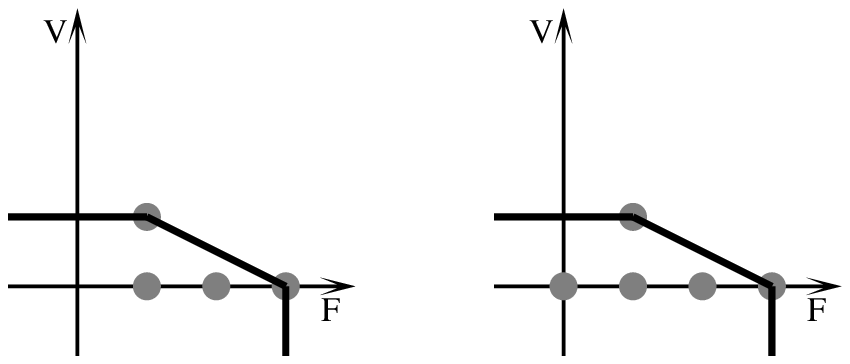}
\end{center}

\noindent From the Newton polygons, $-1/2$ is the candidate of the first slope.
Next, we compute the Gr\"obner basis for the weight $(0,-2,1,3) =
(0,0,1,1)+2(0,-1,0,1)$. The Gr\"obner basis is {\footnotesize
\begin{verbatim}
[    x1*Dx1+3*x2*Dx2+3 , Dx2-Dx1^3 ]
\end{verbatim}
}
\noindent
The radical is generated by
{\footnotesize \begin{verbatim}
[ -x1*Dx1-3*x2*Dx2 , -Dx1^3+Dx2 , 3*x2*Dx1^2*Dx2+x1*Dx2 ,
 -9*x2^2*Dx1*Dx2^2+x1^2*Dx2 , 27*x2^3*Dx2^3+x1^3*Dx2 ]
\end{verbatim}
} \noindent
Here, ${\tt x1} = y_1$, ${\tt x2} = \tau^*$,
${\tt Dx1} = y^*_1$,  ${\tt Dx2} = -\tau$.
It is not bi-homogeneous and then $-1/2$ is a geometric and algebraic slope.

By looking at the two Newton polygons of
\verb@ x1*Dx1+3*x2*Dx2+3@, \verb@ Dx2-Dx1^3@,
we see that there is no more slope that is larger than $-1/2$.
Then, the ACG algorithm terminates here.
\end{example}

The ACG algorithm requires a repetition of Gr\"obner basis
computations in the Weyl algebra of $2n$ variables
to evaluate the slopes of $A_n/I$.
However, if $x_i=0$ is non-micro-characteristic
of type $L=(-r)F+V$ for all $-\infty < r < 0$
and the restriction of
$A_n/I$ to $x_i=0$ is singly generated, we can preprocess the
input so that the input ideal for the ACG algorithm lies in the
Weyl algebra of $2(n-1)$ variables.
The correctness of the following algorithm can be shown by
Laurent and Mebkhout's theorem \ref{Cauchy-Kow}.


\begin{algorithm} \rm (Computing slopes with a preprocessing) \\
\label{algorithm:accelerated-ACG} \\
Step 1: check if $x_i=0$ is non-micro-characteristic of
$A_n/I$ for all types $L$ by calling Algorithm
\ref{algorithm:range}. \\
Step 2: Compute the restriction of $A_n/I$ to
$x_i=0$ and check if it is expressed as $D'/I'$ where $D'$ is the
Weyl algebra of $2(n-1)$ variables. \\
Step 3:
If we failed either in Step 1 or in Step 2,
then apply the ACG algorithm for $I$.  \\
If we succeeded both in step 1 and in step 2,
then try to reduce more variables or apply the ACG algortihm for $I'$.
\end{algorithm}

We can compute restrictions of a given ${\cal D}$-module
by using Oaku's algorithm \cite{O}.
This algorithm is implemented in computer algebraic systems
Macaulay2 and Kan
\cite{Macaulay2}, \cite{openxm}, \cite{Tsai}.
Therefore, the remaining
algorithmic question for the preprocessing is to determine the
range of type $L=(-r)F+V$ for which $x_i=0$ ($ i \leq n-1$) is
non-micro-characteristic.
It follows from the definition of non-micro-characteristic
that the question is nothing but to find the segment
$(-\infty,r_1)$ such that
$$ {\cal V}(\sigma^{(-r)F+V}(I), y_1^*, \ldots,
    y_{i-1}^*, y_i, y_{i+1}^*, \ldots, y_{n-1}^*, \tau^*)
\subseteq {\cal V}(y_1^*, \ldots, y_{n-1}^*, \tau^*) = T^*_\Lambda \Lambda$$
for $r \in (-\infty,r_1)$.
Here, ${\cal V}(f_1, \ldots, f_m)$ is the affine variety defined by
the polynomials $f_1, \ldots, f_m$.
The inclusion condition can algebraically rephrased as
$$ \sqrt{ {\rm in}_{(-r)F+V}(I),
  \xi_1, \ldots, \xi_{i-1}, x_i, \xi_{i+1}, \ldots, \xi_{n} }
 \ni \xi_i.
$$
Since the Gr\"obner fan is a finite union of Gr\"obner cones
\cite{ACG}, the range can be determined by a similar method with
the ACG algorithm.

\begin{algorithm}\label{algorithm:range} \rm \ \\
${\tt range\_of\_nonMC}(H,r_0)$ \\
Input: $H$ is a finite set in $A_n$,
$r_0$ is a negative number or $-\infty$ \\
Output: $r_1$ such that
$x_1=0$ is non-micro-characteristic of type $(-r)F+V$ \\
for $A_n/A_n \cdot \{H\}$, for $ r \in [r_0,r_1)$ . \\
${\tt previous\_slope} = {\tt slope} = r_0$; \\
${\tt F}=(0, \ldots, 0, 1,\ldots, 1)$;
${\tt V} = (0, \ldots, 0, -1,0, \ldots, 0, 1)$; \\
{\tt while} $({\tt slope} != 0)$  {\tt \{} \\
\mbox{\,}\quad ${\tt p} = {\rm numerator}(|{\tt slope}|)$;
               ${\tt q} = {\rm denominator}(|{\tt slope}|)$; \\
\mbox{\,}\quad ${\tt L} = {\tt p F} + {\tt q V}$; \\
\mbox{\,}\quad ${\tt G} = \mbox{ a Gr\"obner basis of }
     H \mbox{ with respect to } <_{L}$; \\
\mbox{\,}\quad
 {\tt if} $\sqrt{ \langle \sigma^{{\tt L}}(H), y_1, y^*_2,
                 \ldots, y^*_{n-1}, \tau^* \rangle } \ni y_1^*$ and \\
\mbox{\,}\quad \ \ $\sqrt{ \langle \sigma^{{\tt pF}+
 ({\tt q}+\varepsilon){\tt V}}(H), y_1, y^*_2, \ldots,
   y^*_{n-1},\tau^* \rangle } \ni y_1^*$ \\
\mbox{\,}\quad {\tt then} {\tt \{} \\
\mbox{\,}\quad \quad   ${\tt previous\_slope} = {\tt slope}$; \\
\mbox{\,}\quad \quad
   ${\tt slope} = {\rm the\ minimum}\ \mbox{ of } \ 0 \ \mbox{ and }$ \\
\mbox{\,}\quad \quad\quad $\{ \mbox{ the slopes $r$
of the Newton polygon } N(P) \,|\, P \in G, r > {\tt previous\_slope} \}$ \\
\mbox{\,}\quad {\tt \}} {\tt else}
{\tt \{} \\ \mbox{\,}\quad\quad return({\tt slope}); \\
\mbox{\,} \quad {\tt \}} \\
{\tt \}} \\
{\tt return}(0); \\
Here, $\varepsilon$ is a sufficiently small positive rational
number so that $pF + (q+\varepsilon) V$ lies in the interior
of a Gr\"obner cone.
\end{algorithm}

\noindent
When $r_0 = -\infty$, we use the convention
${\rm numerator}(|r_0|) = 1$
and
${\rm denominator}(|r_0|) = 0$
and
$F$-non-micro-characteristic means that it is non-characteristic
in the classical sense.

\begin{example} \rm
Suppose $n \geq 2$. ${\tt range\_of\_nonMC}(\{ \pd{1}^2-\pd{2}
\}, -\infty)$ returns $0$.
\end{example}
If the function ${\tt range\_of\_nonMC}(I, -\infty)$ returns $0$,
then $x_1=0$ is non-micro-characteristic of type $p L + q V$ ($p >
0$, $q \geq  0$) for $A_n/I$. We note that it is not always a
clever strategy to call the function with the full set of
generators. In fact, if $x_1=0$ is non-micro-characteristic of
type $L$ for $A_n/J$, then it is non-micro-characteristic of type
$L$ for $A_n/I$ for any $I \supseteq J$. Therefore, for step 1, it
is sometimes more efficient to call the function ${\tt
range\_of\_nonMC}$ for a subset of the generators of the input
ideal as we will see in case of GKZ hypergeometric ideals in the
next section.

\section{Computing slopes of
${\cal H}_{(1,a_2, \cdots , a_n)}(\beta)$}
\label{section:computingSlopes}
\setcounter{equation}{0}

Put
$A=(1,a_2, \cdots ,a_n)$, $a_1=1 < a_2 < \cdots < a_n$.
We will evaluate the slopes of the GKZ hypergeometric
ideal $H_A(\beta)$ associated to the $1 \times n$ matrix $A$
and $\beta \in {\bf C}$ by using the general algorithm given
in Section \ref{section:general-algorithm}.
To apply this algorithm, we need to find non-microcharacteristic varieties
and compute the restrictions of ${\cal H}_A(\beta)$ to these varieties.
When $f_1,\ldots,f_m$ are polynomials in
${\bf C}[x_1,\ldots,x_n,\xi_1,\ldots,\xi_n]$,
we denote by
${\cal V}(f_1,\ldots,f_m)$ the affine subvariety in ${\bf C}^{2n}$
defined by the $f_i$.

The following theorem can be shown by a standard method of Koszul complex
(\cite{Adolphson}, \cite{GKZ}).   

\begin{theorem}
The characteristic variety of ${\cal H}_A(\beta)$ is
${\cal V}(\xi_1, \ldots, \xi_{n-1},x_n \xi_n)$.
In particular, the singular locus of ${\cal H}_A(\beta)$ is $x_n=0$.
\end{theorem}

Note that there is no slope along $x_i=0$, $1 \leq i \leq n-1$,
which can be shown easily.  

Recall that $F=(0,\ldots,0,1,\ldots,1)$ and
$V=(0,\ldots,0,-1,0,\ldots,0,1)$.
For positive number $p$ and
a non-negative number $q$, we define a weight vector $L = p F
+ q V$.
In section \ref{section:microcharacteristic},
we explained the notion of non-micro-characteristic.
When the variety is $y_i=x_i=0$, this notion is rephrased as
follows: for a given left $A_n$-module $A_n/I$, the hyperplane
$y_i=0$ ($1\leq i \leq n-1$) is called
non-micro-characteristic of type $L$ when
$$ \sqrt{\langle \sigma^L(I), y_i,y^*_j, (j\not=i), \tau^* \rangle}
\ni y^*_i. $$

\begin{proposition} \label{proposition:micro}
For the hypergeometric
$A_n$-module ${\cal H}_A(\beta)$,
the variety $y_i =0$, {\rm (}$1 \leq i \leq n-2${\rm )}
is non-micro-characteristic of type $L$ for all
$L = pF + qV$, $p>0$.
\end{proposition}

\begin{proof}
Consider $\partial_i^{a_j}-\partial_j^{a_i} \in H_A(\beta)$.
For all $L$ and
for $i < j \leq n-1$, we have
$\sigma^L(\pd{i}^{a_j}-\pd{j}^{a_i}) = (y^*_i)^{a_j}$, which implies that
$y_i = 0$ is non-micro-characteristic of type $L$.
\end{proof}

Now, let us apply the second step of the algorithm to evaluate
the slopes, i.e., we will compute the restriction of
${\cal H}_A(\beta)$ to $y_i=x_i=0$, $( 1 < i \leq n-2 )$.
Let $s$ be an indeterminate.
Consider the ideal $H_A[s]$ in $A_n[s]$ generated by
$A \theta - s$ and $I_A$.

\begin{theorem}  \label{theorem:restriction}
We have a left $D'[s]$-module isomorphism
\begin{equation}   \label{restriction}
  A_n[s]/(A_n[s] H_A[s]+x_i A_n[s]) \simeq D'[s]/D'[s] H_{A'}[s] , \quad
  i \not= 1
\end{equation}
Here, $D' = {\bf C}\langle x_1, \ldots, x_{i-1}, x_{i+1},
\ldots, x_n, \pd{1}, \ldots, \pd{i-1}, \pd{i+1}, \ldots,
\pd{n} \rangle $ and
$A'=(1,a_2,\ldots,a_{i-1},a_{i+1},\ldots,a_n)$.
\end{theorem}

\begin{proof} \rm
Fix an order $\prec$
such that
$\pd{i} \succ \pd{n} \succ \pd{n-1} \succ \cdots \succ \pd{1}$.
Then,
$$ \pd{n} - \pd{1}^{a_n}, \
   \pd{n-1} - \pd{1}^{a_{n-1}}, \
   \cdots,
   \pd{2} - \pd{1}^{a_{2}}
$$ is the reduced Gr\"obner basis of $I_A$  with respect to
$\succ$.

For each $i$  $(2 \leq i \leq n)$,  by applying the same
method with the proof of \cite[Th 3.1.3]{SST}, we can prove
that ${\rm in}_{(-e_i,e_i)}(H_A(s))$ is generated by
$A \theta - s$ and
${\rm in}_{e_i}(\pd{j} - \pd{1}^{a_j})$, $(j=2, \ldots,n)$.

Define the $V_i$-filtration $F_k[s]$ of $A_n[s]$ by
$$ F_k[s] = \left\{ \sum a_{\alpha\beta\gamma} x^\alpha \pd{}^\beta s^\gamma
\,|\, (\beta-\alpha) \cdot e_i \leq k \right\}.$$
We compute the
restriction of $A_n[s]/H_A(s)$ to $x_i=0$ by Oaku's algorithm
(see, e.g., \cite[Theorem 5.2.6, Alg. 5.2.8]{SST}). Since
${\rm in}_{(-e_i,e_i)}(x_i \pd{i} - x_i \pd{1}^{a_i}) =
  x_i \pd{i}$,
the $b$-function of $H_A(s)$ along $x_i=0$ is $b(p)=p$.
Therefore, by \cite[Theorem 5.2.6]{SST},
we have the following isomorphism as left $D'[s]$-modules:
\begin{eqnarray*}
& &A_n[s]/(H_A(s)+x_i A_n[s])  \\
&\simeq&
F_0[s]/(F_0[s](A\theta-s)
 + \sum_{j=2, j\not=i}^n F_{0}[s](\pd{j}-\pd{1}^{a_j})
 + F_{-1}[s](\pd{i}-\pd{1}^{a_i})+x_i F_1[s]) \\
&\simeq& D'[s]/\left(D'[s](A\theta-s) +\sum_{j=2, j\not=i}^n
         D'[s](\pd{j}-\pd{1}^{a_j}) \right)
\end{eqnarray*}
\end{proof}

We can specialize $s$ to any complex number $\beta$.
So, we have
$$ A_n/(A_nH_A(\beta)+x_i A_n) \simeq
D'/D'H_{A'}(\beta),$$
which means that the restriction of ${\cal H}_A(\beta)$ to $x_i = 0$
can be exactly expressed in terms of the GKZ system for
smaller $A$.
By applying our algorithm \ref{algorithm:accelerated-ACG}
of computing slopes by
reduction of number of variables to variables
$x_2, \ldots, x_{n-2}$,
we obtain the following theorem
by Proposition \ref{proposition:micro}
and Theorem \ref{theorem:restriction}.
\begin{theorem}
The  geometric slopes of ${\cal H}_A(\beta)$ along $x_n=0$ at the
origin and ${\cal H}_{(1,a_{n-1},a_n)}(\beta)$ along $x_3=0$ at the
origin coincide.
\end{theorem}

\begin{example} \rm
We note that
$\mbox{(the algebraic slopes)} \not=
\mbox{ (the geometric slopes)}$
in general.
For example,
let us apply the ACG algorithm to get the algebraic slopes
of $H_A(-30)$ for $A=(1, 3, 7)$.
This ideal is generated by
$$x_1\pd{1}+3 x_2\pd{2}+7x_3\pd{3}+30,
 \pd{1}^3-\pd{2} , -\pd{1}^2\pd{2}+\pd{3} , \pd{2}^3-\pd{1}^2\pd{3}.
$$ The output is $$\{ -1,  -3/4,  -1/2 \}.$$
On the other hand, if we apply the ACG algorithm to get the
geometric slopes, the output is $\{ -3/4 \}$.
\end{example}

\begin{example} \rm
$(\mbox{the slopes of } {\cal H}_{(1,a_{n-1},a_n)}(\beta)) \not=
(\mbox{the slopes of } {\cal H}_{(1,a_n)}(\beta))$   in general. \\
Let us take the example: $A = (1,3,7)$.
Consider the
hypergeometric $A_3$-module $A_3/I$ where $I = H_{(1,3,7)}(-30)$.
As we have seen, the slope of this system along $x_3=0$
is $\{ -3/4 \}$.

Consider ${\rm in}_L(I)$ for $L =  F + 4 V$.
By computing Gr\"obner basis with respect to $L$,
we can see that
$$ {\cal V}({\rm in}_L(I)) = {\cal V}(x_2, \xi_1, \xi_3) \cup
                      {\cal V}(\xi_1, \xi_2, \xi_3).
$$
It is not included in ${\cal V}(\xi_2)$
Hence $x_2=0$ is micro-characteristic of type $L$
and we cannot apply for the restriction criterion.

The condition ``non-microcharacteristic
for all the filtration $pF+qV$'' cannot be taken a way to
evaluate the slopes by the restriction.
In fact, it can be
easily checked by the ACG algorithm that the set of the
geometric slopes of ${\cal H}_{(1,a_n)}(\beta)$ is equal to
$\{ 1/(1-a_n) \}$.
Hence, the set of the slopes of
${\cal H}_{(1,7)}(\beta)$ is $\{ -1/6 \}$, which is not equal to
$\{ -3/4 \}$.
\end{example}

We have shown that the computation of the slopes of
${\cal H}_A(\beta)$ is reduced to the three variables case.
The slopes of this case are as follows.

\begin{theorem}
$$(\mbox{\rm the slopes of } {\cal H}_{(1,a_{n-1},a_n)}(\beta))
= \{ a_{n-1}/(a_{n-1}-a_{n}) \}$$
\end{theorem}

\begin{proof} \rm
We fix  some notation:
\begin{enumerate}
\item $A=(1, a , b)\in {\bf Z}^3$ and $1 < a < b$.
\item $P_1=\partial_1^a - \partial_2$, $P_2=\partial_1^b-\partial_3$,
$P_3=\partial_2^b - \partial_3^a$, $P_4=x_1\partial_1 + a
x_2\partial_2 + b x_3\partial_3 - \beta$
\item Let $\Lambda$ be the linear form with slope $-a/(b-a)$
(i.e. $L = a F + (b-a) V$).
\item Let $L,L'$ be  linear forms. We say that $L > L'$
if ${\rm slope}(L) > {\rm slope}(L')$.
\item We will write  $y_1=x_1$, $y_2=x_2$, $t=x_3$.
\end{enumerate}

The operators $P_1,P_2,P_3,P_4$ are in $H=H_A(\beta)$.
Let $\Lambda$ be $a_{n-1}/(a_n-a_{n-1})$.
Then, we have the following claims.
\begin{enumerate}
\item For all linear form $L$ we have $\sigma^L(P_1)=(\eta_1^*)^a$
and so $\eta_1^* \in \sqrt{\sigma^L(H)}$ for all $L$.
\item For all
linear form $L$ we have $\sigma^L(P_4) = y_1 \eta_1^* + a y_2
\eta_2^* + b \tau^* (-\tau)$
\item For all linear form $L > \Lambda$ we have $\sigma^L(P_3) =
(\tau)^a$ and so $\tau \in \sqrt{\sigma^L(H)}$ for all $L> \Lambda$.
\item So, for all $L> \Lambda$ we have ${\rm Ch}^L({\cH}) \subset
T^*_{y_2=0}\CC^3 \cup T^*_{\CC^3}\CC^3$ and then
$\sqrt{\sigma^L(H)}$ is bi-homogeneous and $L$ is not a geometric
slope of $\cH$.
\item On the other hand we have, for $L<\Lambda$,
$\sigma^L(P_3) = (\eta_2^*)^b$. Then $\eta_2^* \in \sqrt{\sigma^L(H)}$
and ${\rm Ch}^L(\cH) \subset T^*_{t=0}\CC^3 \cup T^*_{\CC^3} \CC^3$. So
$L$ is not a geometric slope of $\cH$.
\item So, the only possible geometric slope of $\cH$ is $\Lambda$.
\end{enumerate}

Now, suppose that $\Lambda$ is not a slope.
Then, there is no slope,
which implies that the $L$-characteristic variety
${\rm Ch}^L({\cal H}_A(\beta))$
is the same for all
$L = p F + q V$, $p, q > 0$
by \cite{ACGjpaa} and \cite[Th 3.4.1]{Lau-ens}. 
It follows from
$$[P_1, P_2] = 0, \ [P_1,P_4] = a P_1,\  [P_2, P_4] = b P_2 $$
and the Buchberger algorithm
that
$\{ P_1, P_2, P_4\}$ is a Gr\"obner basis for the order
defined by the weight vector
$L = (0,0,-N,1,1,N+1)$, $N \geq b$
and a tie-breaking term order such that
$ x_2 \succ \pd{3} \succ \pd{1} \succ \pd{2}$.
Therefore, the initial ideal ${\rm in}_L(H_A(\beta))$ is generated by
$ \xi_1^a$, $\xi_3$, $x_1 \xi_1 + a x_2 \xi_2 + b x_3 \xi_3$
and hence the $L$-characteristic variety is equal to
$T^*_{y_2=0}\CC^3 \cup T^*_{\CC^3}\CC^3$.
This fact contradicts that the $L$-characteristic variety is the same for
all $L$.
\end{proof}
\medbreak


Finally, let us remark on an analytic meaning of slopes.
Let $X={\bf C}^n$ and $Y=\{x\vert x_n=0\} \subset X$.
We denote by
${\cal O}_{\widehat{X\vert Y}}
 = {\bf C}\{ x_1, \ldots, x_{n-1} \} [[x_n]]$
the formal completion of
${\cal O}_X$ along $Y$.
For each real number $s\in [1,+\infty)$
we denote by ${\cal O}_{{X\vert Y}}(s)$ the
subsheaf of ${\cal O}_{\widehat{X\vert Y}}$ of Gevrey functions of order
$s$
(along $Y$). The sheaf ${\cal O}_{{X\vert Y}}(1)$ is the restriction
${\cal O}_{{X\vert Y}}$
and, by definition, we write
${\cal O}_{{X\vert Y}}(+\infty) = {\cal O}_{\widehat{X\vert Y}}$.
For any holonomic ${\cal D}_X$-module $M$, Mebkhout associates in
\cite{Meb-positivite} the sheaf ${\rm Irr}_Y(s)(M)$ as the solution sheaf
${\bf R}Hom_{\cal D}(M,{\cal O}_{{X\vert Y}}(s)/{\cal O}_{{X\vert Y}})$.
One fundamental result in the irregularity of $\cal D$-modules is the fact
that
${\rm Irr}_Y(s)(M)$ is a perverse sheaf, for any $s$ (see
\cite{Meb-positivite}).
These sheaves define a filtration of the irregularity of $M$ along $Y$,
i.e.
${\rm Irr}_Y(M):={\rm Irr}_Y(+\infty)(M)$.
The main result of \cite{Lau-Meb-papa} is
that
$1/(1-s)$ is a slope of $M$ w.r.t. $Y$ if and only if $s$ is
a gap of the graduation defined
by the filtration on the irregularity.
In other words, $1/(1-s)$ is a slope if and only if
${\rm Irr}_Y(s)(M)/{\rm Irr}_Y(<s)(M) \not= 0$.


\section{Rational Solutions and Reducibility}
\setcounter{equation}{0}

Our ultimate aim of studying the slopes of $H_A(\beta)$
is to get a better understanding on solutions of this system.
We are far from the goal, but
to this end, it will be useful to present some facts
on classical solutions and a relation to generalized confluent hypergeometric
functions.

In case of the hypergeometric ideal
associated to homogeneous monomial curves,
Cattani, D'Andrea, Dickenstein \cite{CDD}  studied
rational solutions and reducibility of the system.
We will study rational solutions and reducibility of our system.

\begin{theorem}
Any rational solution of the hypergeometric system
$H_{(1,a_2, \ldots, a_n)}(\beta)$ is a polynomial.
It has a polynomial solution if and only if
$\beta \in {\bf N}=\{0, 1, \ldots, \}$.
The polynomial solution is the residue of
$\exp\left(\sum x_i t^{a_i}\right) t^{-\beta}$
at the origin $t=0$:
$$  \int_C \exp\left(\sum x_i t^{a_i}\right) t^{-\beta} \frac{dt}{t}. $$
Here, $C$ is a circle that encircles the origin in the positive direction.
\end{theorem}

\begin{proof}
Since the singular locus of
$H_A(\beta)$,
($A=(a_1,a_2, \ldots, a_n)$, $a_1=1$) is
$x_n=0$,
any rational solution $f$ is a Laurent polynomial
with poles on $x_n=0$.
Take a weight vector
$w = (0, 1,1, \ldots, 1)$.
Then, we have
${\rm in}_w(I_A) = \langle \pd{2}, \ldots, \pd{n} \rangle $.
The initial term ${\rm in}_w(f)$ is annihilated by
$$ \sum a_i \theta_i - \beta, \quad \pd{2}, \ldots, \pd{n}. $$
Therefore,
$x_1^\beta = {\rm in}_w(f)$.
This implies $\beta \in {\bf N}$, because $f$ has a pole only on
$x_n=0$.

Take a ${\bf Z}$-basis of
${\rm Ker}\, ({\bf Z}^n \stackrel{A}{\rightarrow} {\bf Z})$
as
$$(-a_2, 1,0, \ldots, 0), (-a_3, 0,1,0, \ldots, 0), \ldots,
 (-a_n, 0, \ldots, 0, 1)$$
to construct series solutions.
Since \cite[Prop. 3.4.1]{SST} holds for non-homogeneous $A$ as well,
formal series solution $g$ of $H_A(\beta)$ satisfying
${\rm in}_w(g) = x_1^\beta$
can be uniquely expressed as
\begin{equation} \label{formal_sol}
\sum_{m \in{\bf N}^{n-1}}
  \frac{\beta (\beta-1) \cdots (\beta-\sum m_k a_k + 1)}
       {m_2 ! \cdots m_n !}
  \left( \frac{x_2}{x_1^{a_2}} \right)^{m_2} \cdots
  \left( \frac{x_n}{x_1^{a_n}} \right)^{m_n} x_1^\beta.
\end{equation}
When $\beta \in {\bf N}$, it is a polynomial.
The rest part of the theorem is easy to show.
\end{proof}

\bigbreak

Let $R$ be the ring of differential operators of
$n$ variables with rational function coefficients over ${\bf k}={\bf C}$.
A left ideal $J$ of $R$ is called
{\it irreducible} when $J$ is a maximal ideal in $R$.
We will study the reducibility of $R\cdot H_A(\beta)$.

We assume that $J$ is zero-dimensional, i.e., $r = {\rm dim}_{{\bf
C}(x)}\, R/J < +\infty$. Let $V=V(J)$ be the vector space of
holomorphic solutions of $J$ on a simply connected open set
contained in the non-sigular domain of $J$. It is known that ${\rm
dim}_{\bf C}\, V = r$. Define $I(V)$ by $R\cdot\{ \ell \in R\,|\,
\ell \bullet f = 0 \ \mbox{ for all }\ f \in V \}$. If $J \subset
I(V)$, $J \not= I(V)$, then we have $ {\rm dim}_{{\bf C}(x)} R/J <
{\rm dim}_{{\bf C}} V = r$ because of the zero-dimensionality of
$J$. Therefore, we have
$$ J = I(V(J)). $$
Under this
correspondence of ideals and solutions, a zero-dimensional ideal
$J$ of $R$ is reducible if and only if there exists a proper
subspace $W$ of the solution space of $V(J)$ such that
$ 0 < {\rm dim}_{{\bf C}(x)} R/I(W) < {\rm dim}_{{\bf C}(x)} R/J$.
In case
of one variable, the reducibility is equivalent to saying that the
generator of the ideal can be factored in $R$.

\bigbreak
\begin{theorem}
The systems of differential equations $R \cdot H_A(\beta)$
is reducible if and only if $\beta \in {\bf Z}$.
\end{theorem}

\begin{proof}
Any curve is Cohen-Macaulay.
By applying the theorem of Adolphson \cite{Adolphson},
the holonomic rank of $H_A(\beta)$ is $a_n$ for all $\beta$.

Put $M(\beta) = A_n/H_A(\beta)$. Consider the left $A_n$-morphism
\begin{equation} \label{eq:beta-isom}
 \pd{1}\ :\ M(\beta) \rightarrow M(\beta+1).
\end{equation}
It has the inverse when $\beta \not= -1$.
Therefore,  we have
$M(-1) \simeq M(-2) \simeq M(-3) \simeq \cdots$
and
$M(0) \simeq M(1) \simeq M(2) \simeq \cdots$.

When $\beta \in {\bf N}$, the system
admits polynomial solution,
then it is reducible.
It is also easy to see that when $\beta \in {\bf Z}_{<0}$,
the equation is reducible.
In fact, consider the left $D$-morphism
$$ \pd{1}\ :\ M(-1) \rightarrow M(0). $$
It induces a morphism to the solutions
by
$$ f \rightarrow \pd{1} \bullet f. $$
The solution $f=1$ of $M(0)$ is sent to zero,
so the image of $\pd{1}$ gives a proper subspace of solutions
in the solution space of $M(-1)$.
To find differential equations for the subspace,
take all $\ell$ such that $\ell \pd{1} \in H_A(0)$.
Then, $\{ \ell \} \subset H_A(-1)$.
By the isomorphism (\ref{eq:beta-isom}),
we conclude that when  $\beta \in {\bf Z}_{<0}$, the system
is reducible.

Let us prove that the system is irreducible  when $\beta \not\in {\bf Z}$
by applying the result of Beukers, Brownawell, and Heckman \cite{BBH}.
For this purpose, we firstly construct convergent
series solution of $H_A(\beta)$.
Take $w = (1,1, \ldots, 1,0)$.
Then the degree of ${\rm in}_w(I_A)$ is equal to $a_n$.
Since $I_A$ contains the elements of the form
$\underline{\pd{i}^{a_n}} - \pd{n}^{a_i}$,
the radical of
${\rm in}_w(I_A)$ is $ \langle \pd{1}, \ldots, \pd{n-1} \rangle$.
Therefore, the top dimensional standard pairs have the form
$(\pd{}^b, \{n\})$.

Let $v$ be the zero of the indicial ideal associated to
$(\pd{}^b, \{n\})$:
$$ v_1 = b_1, \ldots, v_{n-1} = b_{n-1},
   v_n = \frac{\beta-\sum a_i b_i}{a_n}.
$$
Assume $\beta \not\in {\bf Z}$ or $\beta \gg 0$.
Taking the lattice basis
$(a_2, -1, 0, \ldots, 0)$, \ldots, $(a_n, 0, \ldots, 0, -1)$,
we have the following $a_n$ linearly independent
convergent series solutions
\begin{eqnarray}
& &
\sum_{m \in {\bf N}^{n-1}}
 \frac{ v_2 (v_2-1) \cdots (v_2-m_2+1) \cdots  v_n (v_n-1) \cdots (v_n-m_n+1)}
      { (v_1+1) (v_1+2) \cdots (v_1 + \sum a_k m_k)} \nonumber \\
& &\quad\quad\quad\quad\quad\quad\quad\quad\cdot
\left( \frac{x_1^{a_2}}{x_2} \right)^{m_2} \cdots
 \left( \frac{x_1^{a_n}}{x_n} \right)^{m_n} x^v.  \label{convergent_sol}
\end{eqnarray}

We consider the change of variables: $$ y_1 = x_1, y_2 =
x_1^{a_2}/x_2, \ldots, y_n = x_1^{a_n}/x_n. $$ The inverse of this
change of variables is also rational and the change of variables
induces that of $\pd{i}$. We denote by $\Phi$ the operation of
these change of variables of $x_i$ and $\pd{i}$. Since the
irreducibility is invariant under any birational change of
variables, we will prove the irreduciblity of the ideal $J = R
\cdot \Phi(x^{-v} H_A(\beta) x^v)$ where $R = {\bf C}(y) \langle
\pd{y_1}, \ldots, \pd{y_n} \rangle$.

Let $V$ be the solution space of $J$ spanned by the series
$\Phi(\mbox{(\ref{convergent_sol})}\cdot x^{-v})$ near $y=0$. If
$J=I(V)$ is reducible, then there exists a proper subspace $W$ of
$J$ such that $0 < {\rm dim}_{{\bf C}(y)}\, R/I(W) < a_n$. We
consider the vector space $W' = \{ \, f(0, \ldots, 0, y_n) \,|\, f
\in W \}$. It is easy to see that $ {\rm dim}_{{\bf C}(y_n)}
R'/I(W') \leq {\rm dim}_{{\bf C}} W $ where $R' = {\bf C}(y_n)
\langle \pd{y_n} \rangle$. Let us prove that ${\rm dim}_{{\bf C}}
W' = a_n$ when $\beta \not\in {\bf Z}$, which implies the
irreduciblity of $J$ by a contradiction.

We restrict the series $\Phi(\mbox{(\ref{convergent_sol})} \cdot x^{-v})$ to
$y_1=x_1=0, y_2=x_1^{a_2}/x_2 = 0$, \ldots, $y_{n-1}=x_1^{a_{n-1}}/x_n = 0$
and replace $y_n=x_1^{a_n}/x_n$ by $z$.
Without a loss of generality, we may assume $v_1 =0$.
Then the restricted series has the form
\begin{equation} \label{eq:convergent:one-variable}
\sum_{m=0}^\infty
 \frac{ v_n (v_n-1) \cdots (v_n-m+1)}
      { (a_n m) ! }(-z)^m.
\end{equation}
It is annihilated by the ordinary differential operator
$$ (a_n \theta_z) (a_n \theta_z -1) \cdots (a_n \theta_z-a_n + 1) -
   z (\theta_z + v_n).$$
By replacing $z/a_n^{a_n}$ by $x$, we obtain the generalized
hypergeometric ordinary differential equation
\begin{equation} \label{eq:generalizedHG}
 \theta_x (\theta_x - 1/a_n) \cdots (\theta_x - (a_n-1)/a_n)-
x (\theta_x + v_n).
\end{equation}
By \cite{BBH}, this ordinary differerential equation of rank $a_n$
is reducible if and only if $v_n - k/a_n \not\in {\bf Z}$ for all
$k=0,1, \ldots, a_n-1$. If one of them is an integer, $\beta$
becomes an integer. Therefore, the ideal $I(W')$ contains the
principal ideal generated by (\ref{eq:generalizedHG}), which is
maximal when $\beta \not\in {\bf Z}$. We conclude that $I(W')$ is
generated by (\ref{eq:generalizedHG}) and hence ${\rm dim}_{\bf
C}\, W' = a_n$.
\end{proof}

\bigbreak
\noindent
{\it Acknowledgement\/}:
The main part of this paper was obtained from discussions during
second author's visit to Universidad de Sevilla in July, 2000.
The authors are grateful to the university of providing us the
opportunity of the intensive discussions.

\end{document}